\documentclass{amsart}

\usepackage{amsfonts}
\usepackage{amsmath}
\usepackage{amsthm}
\usepackage{amssymb}
\usepackage[english]{babel}
\usepackage{graphics}
\usepackage{latexsym}
\usepackage{longtable} \usepackage{mathrsfs}
\usepackage{graphicx}
\usepackage{wrapfig} 
\usepackage[pdftex,colorlinks=true]{hyperref}

\newtheorem{theorem}{Theorem}

\newtheorem{lemma}[theorem]{Lemma}
\newtheorem{proposition}[theorem]{Proposition}

\newtheorem{claim}[theorem]{Claim}
\newtheorem{example}[theorem]{Example}

\theoremstyle{definition}
\newtheorem{definition}[theorem]{Definition}
\newtheorem{remark}[theorem]{Remark}

\newcommand{\mH}{\mathcal{H}}

\newcommand{\A}{\textrm{A}}\newcommand{\R}{\mathbb{R}}

\newcommand{\mS}{\mathbb{S}}
\newcommand{\mB}{\mathbb{B}}
\newcommand{\X}{\textbf{X}}

\renewcommand{\P}{\textrm{P}}
\newcommand{\Q}{\textrm{Q}}

\newcommand{\noi}{\noindent}
\newcommand{\ms}{\medskip}

\newcommand{\al}{\alpha}
\newcommand{\be}{\beta}

\newcommand{\de}{\delta}
\newcommand{\De}{\Delta}
\newcommand{\e}{\varepsilon}
\newcommand{\si}{\sigma}

\newcommand{\Om}{\Omega}

\newcommand{\larrow}{\longrightarrow}
\newcommand{\ot}{\otimes}

\newcommand{\ri}{\rightarrow}

\newcommand{\p}{\partial}
\newcommand{\sub}{\subseteq}
\newcommand{\set}{\setminus}
\newcommand{\by}{\times}
\newcommand{\rk}{\textrm{rk}}

\newcommand{\sgn}{\textrm{sgn}}
\newcommand{\ess}{\textrm{ess}}

\newcommand{\Div}{\textrm{Div}}

\newcommand{\spn}{\textrm{span}}

\newcommand{\bt}{\begin{theorem}}\newcommand{\et}{\end{theorem}}
\newcommand{\bd}{\begin{definition}}\newcommand{\ed}{\end{definition}}
\newcommand{\bl}{\begin{lemma}}\newcommand{\el}{\end{lemma}}
\newcommand{\beq}{\begin{equation}}\newcommand{\eeq}{\end{equation}}
\newcommand{\bc}{\begin{claim}}\newcommand{\ec}{\end{claim}}
\newcommand{\bp}{\begin{proof}}\newcommand{\ep}{\end{proof}}

\newcommand{\BPP}{\medskip \noindent \textbf{Proof of Proposition} }
\newcommand{\BPT}{\medskip \noindent \textbf{Proof of Theorem} }

\numberwithin{equation}{section}
\numberwithin{theorem}{section}

\begin{document}

\title{$L^\infty$ Variational Problems for Maps and the Aronsson PDE System}

\author{\textsl{Nikolaos I. Katzourakis}}
\address{BCAM - Basque Center for Applied Mathematics, Mazarredo 14, 48009, Bilbao, Spain}
\email{nkatzourakis@bcamath.org}


\date{}



\begin{abstract} By employing Aronsson's Absolute Minimizers of $L^\infty$ functionals, we prove that Absolutely Minimizing Maps $u:\R^n \larrow \R^N$ solve a ``tangential'' Aronsson PDE system. By following Sheffield-Smart \cite{SS}, we derive $\De_\infty$ with respect to the dual operator norm and show that such maps miss information along a hyperplane when compared to Tight Maps. We recover the lost term which causes non-uniqueness and derive the complete Aronsson system which has \emph{discontinuous coefficients}. In particular, the Euclidean $\infty$-Laplacian is $\De_\infty u = Du \ot Du : D^2u\, +\, |Du|^2[Du]^\bot \De u$ where $[Du]^\bot$ is the projection on the null space of $Du^\top$. We exibit $C^\infty$ solutions having interfaces along which the rank of their gradient is discontinuous and propose a modification with $C^0$ coefficients which admits \emph{varifold solutions}. Away from the interfaces, Aronsson Maps satisfy a structural property of local splitting to $2$ phases, an horizontal and a vertical; horizontally they possess gradient flows similar to the scalar case and vertically solve a linear system coupled by a scalar Hamilton Jacobi PDE. We also construct singular $\infty$-Harmonic local $C^1$ diffeomorphisms and singular Aronsson Maps. 
\end{abstract}

\maketitle

\section{Introduction} \label{section1}

Let $H \in C^1((\R^N \ot \R^n) \by \R^N \by \R^n)$ be a non-negative Hamiltonian. In this paper we are interested in the theory of \emph{Aronsson PDE systems} which arises when one considers $L^\infty$ variational problems  of the \emph{supremal functional} 
 \beq \label{eq1}
E_{\infty}(u,\Om)\ :\ = \ \underset{x\in\Om}{\ess \, \sup}\ H
\big(Du(x),u(x),(x)\big), \ \ \Om  \sub \R^n,
 \eeq
defined on locally Lipschitz maps $u :\R^n \larrow \R^N$. The study of such problems has been initiated by Aronsson (\cite{A1} - \cite{A7}). The scalar case of $N=1$ has been well studied ever since. By introducing the appropriate minimality notion for \eqref{eq1}, that of \emph{Absolute Minimizers} (recalled in Section \ref{section2}), Aronsson studied solutions $u\in C^2(\R^n)$ of what we now call Aronsson's PDE
 \beq  \label{eq2}
\textrm{A}_\infty u\ := \ H_\P(Du,u,\_)^\top D\big(H(Du,u,\_)\big)\ = \ 0.
 \eeq
Here $H=H(\P,\eta,x)$, subscripts of $H$ denote derivatives and vectors are viewed as columns. \eqref{eq2} arises as an ``Euler-Lagrange'' PDE governing Absolute Minimizers of \eqref{eq1}. The technology of Viscosity Solutions allowed the rigorous study of the generally singular solutions (see e.g. \cite{K}) and a fairly complete scalar theory has been developed (see e.g.\ \cite{J, BJW1, BJW2, CEG, CWY, Y, ACJS, C} and references therein). In particular, for $N=1$ and $H(p)=\frac{1}{2}|p|^2$, there is a triple equivalence among viscosity solutions $u\in C^{0,1}(\R^n)$ of the $\infty$-Laplacian
 \beq  \label{eq3}
 \De_\infty u\ := \  Du \ot Du :D^2 u\ \equiv\ D_i u\, D_j u\, D^2_{ij}u \ = \ 0,
 \eeq
Absolute Minimizers of $E_\infty(u,\Om)=\frac{1}{2}\|Du\|^2_{L^\infty(\Om)}$ and the so-called \emph{Optimal Lipschitz Extensions}, function $u\in C^{0,1}(\R^n)$ satisfying $\textrm{Lip}(u,\Om)=\textrm{Lip}(u,\p \Om)$ for all $\Om \subset \subset \R^n$, where $\textrm{Lip}$ is the Lipschitz functional
 \beq \label{eq4}
\textrm{Lip}(u,K) \ =\ \underset{x,y \in K ,x \neq y}{\sup}
\frac{|u(x)-u(y)|}{|x-y|},\ \ \ K\sub \R^n.
 \eeq
In this work we consider the full vector case of \eqref{eq1} and the related Aronsson PDE systems which arise. The case of general maps has hardly been studied. Recently, Sheffield and Smart \cite{SS} distilled the appropriate notion of optimality for Lipschitz Extensions and strengthened Absolute Minimality to \emph{Tightness}, establishing the equivalence between regular \emph{Tight Maps} $u:\Om \sub \R^n \larrow \R^N$ satisfying 
 \beq
\sup_{\{\textrm{L}u>\textrm{L}v\}}\textrm{L}u\ \leq\ \sup _{\{\textrm{L}v>\textrm{L}u\}}\textrm{L}v, \ \ \ v=u \text{ on }\p \Om,
\eeq
(where $\textrm{L}u(x):=\lim_{r\ri 0}\textrm{Lip}(u,\mB_r(x))$ is the local Lipschitz constant) and solutions of the following \emph{$\infty$-Laplace PDE system}
\beq \label{eq5}
\De_{\infty_T} u\ :=\ D \big(D u\, a\big)a\ = \ D_i \big(D_j u\, a_j \big)\, a_i \ = \ 0.
\eeq
Here $a : \R^n \larrow \mS^{n-1}$ is the unit vector field realizing the \emph{operator norm} on $\R^N \ot \R^n$
\beq \label{eq6}
\|\P\|\ := \ \max_{w\in \mS^{n-1}}|\P w| \ =\ \max \sigma (\P^\top \P)^{\frac{1}{2}} \ = \ \max_{|w|=|\xi|=1} \P : \xi \ot w
\eeq
of the gradient matrix $Du$, that is, $\|Du\|^2=(Du^\top Du): a \ot a$. The usage of the nonsmooth operator norm  on the space of gradients is necessary for the coincidence
\beq \label{eq7}
\underset{\Om}{\ess \,\sup}\, \|Du\| \ = \ \textrm{Lip} (u,\Om)
\eeq 
when $u \in C^1(\R^n)^N$, in order to connect Lipschitz Extensions with $L^\infty$ Variational Problems. In the vector case $N>1$ these two fields decouple, since \eqref{eq7} fails for the Euclidean norm $|\P|=(\P_{\al i}\P_{\al i})^{\frac{1}{2}} = (\P :\P)^{\frac{1}{2}}$ on $\R^N \ot \R^n$. In another recent paper, Capogna and Raich \cite{CR} used the Hamiltonian $H(\P):= {|\P|^n}/{\textrm{det}(\P)}$ in \eqref{eq1} and launched an $L^\infty$ variational approach to extremal \emph{Quasi-Conformal mappings} $u : \R^n \larrow\R^n$. Also, Ou, Troutman and Wilhelm \cite{OTW} and Wang and Ou \cite{WO} studied geometric aspects of ``tangentially'' (as we show herein) $\infty$-Harmonic maps. 

In Section \ref{section2} we employ Aronsson's classical minimality notion for \eqref{eq1} and prove that regular Absolutely Minimizing Maps $u:\R^n \larrow \R^N$ solve the following \emph{Tangential Aronsson PDE system}:
 \begin{align}  \label{eq8}
\big(\textrm{A}_{\top , \infty} u\big)_\al \ :=\ & \Big(H_{\P_{\al i}}(Du,u,\_)\, H_{\P_{\be j}}(Du,u,\_)\Big)\, D^2_{ij}u_\be \nonumber\\
+\ & H_{\P_{\al i}}(Du,u,\_) \, \Big( H_{\eta_\be}(Du,u,\_)D_i
u_\be + H_{x_i}(Du,u,\_)\Big) \ =\ 0.
 \end{align}
The choice $H(\P,\eta,x)= \frac{1}{2}|\P|^2$ in \eqref{eq8} gives a tangential Euclidean version of the $\infty$-Laplacian studied in \cite{OTW}, \cite{WO} which reads 
\beq  \label{eq9}
\big(\De_{\top , \infty} u\big)_\al\ :=\ D_i u_\al \, D_j u_\be \, D^2_{\ij}u_\be \ =  \ 0.
 \eeq
We also characterize regular tangentially $\infty$-Harmonic Maps by means of \emph{gradient flows with parameters}. This result extends directly to solutions of \eqref{eq8}. 

In Section \ref{section3} we follow \cite{SS} and choose $H(\P)=\frac{1}{2}\| \P \|_{*}^2$ as our Hamiltonian, where $\| \_ \|_{*}$ is the dual expression of the operator norm on $\R^N \ot \R^n$
\beq \label{eq10}
\|\P\|_{*}\ := \ \max_{\xi\in \mS^{N-1}}|\xi^\top \P| \ =\ \max \sigma (\P \P^\top)^{\frac{1}{2}} \ = \ \max_{|w|=|\xi|=1} \P : \xi \ot w.
\eeq
By considering the $L^p$ functionals
\beq \label{eq11}
E_p(u,\Om)\ := \ \int_\Om H(Du)^p,
\eeq
we derive in the limit as $p \ri \infty$ the \emph{dual $\infty$-Laplacian which governs Tight Maps}:
\beq \label{eq12}
\De_{\infty_{T^*}} u \ := \ (e^\top Du)\ot (e^\top Du):(e^\top D^2u)\, e\ \, +\, \ e^\bot \Div (e \ot e\, Du) \ =\ 0. 
\eeq
Here $e : \R^n \larrow \mS^{N-1}$ realizes the dual norm $\|Du\|_{*}^2=(Du \, Du^\top): e \ot e$ and $e^\bot :=I-e\ot e$. The reduction of \eqref{eq12} to \eqref{eq3} when $N=1$ is immediate: $e\equiv \pm 1$ and $e^\bot \equiv 0$. The same choice of $H$ in \eqref{eq8} reveals that \emph{Absolutely Minimizing Maps with respect to the (dual) operator norm miss all the information of \eqref{eq12} along the hyperplane $(\spn[e])^\bot \sub \R^N$}.

In Section \ref{section4} we recover the information that \eqref{eq8} fails to encapsulate. For simplicity we take $H=H(Du)$, $H\in C^2(\R^N \ot \R^n)$ and we derive the \emph{complete Aronsson PDE system} (see therein for notation)
\beq \label{eq13}
\A_\infty u\ :=\ \Big(H_\P(Du) \ot H_\P(Du)\ +\ H(Du)[H_\P(Du)]^\bot H_{\P \P}(Du) \Big):D^2u\ = \ 0.
\eeq
Here $[H_\P(Du(x))]^\bot$ is the projection on the null space of $H_\P(Du(x))^\top :\R^N \larrow \R^n$. In particular, for $H(\P)=\frac{1}{2}|\P|^2$, $\De_\infty$ with respect to the Euclidean (Frobenious) norm reduces to
\beq \label{eq14}
\De_\infty u\ :=\ Du\ot Du : D^2u\ +\ |Du|^2[Du]^\bot \De u\ = \ 0.
\eeq
It is quite important however in view of the results in \cite{SS,CR} to work in the generality of \eqref{eq13}, since the most ``suitable" choice of Hamiltonian is not completely settled. The full system \eqref{eq13} remedies the non-uniqueness defects associated to \eqref{eq8} due to the omitted term: every unit speed curve $u:\R \larrow \R^N$ is tangentially $\infty$-Harmonic but not \emph{$\infty$-Harmonic}, since for \eqref{eq14} and $n=1$ solutions are straight lines. However, the full PDE \eqref{eq13} has \emph{discontinuous coefficients even for $C^\infty$ solutions}. There exist smooth $\infty$-Harmonic Maps $\R^2\larrow \R^2$ which are not submersions and their rank is not everywhere equal to $2$ (Subsection \ref{subsection2}), hence generally  they do \emph{not} have constant rank. We therefore propose the following natural modification of \eqref{eq13}:
\beq \label{eq13a}
\Big(H_\P(Du) \ot H_\P(Du)\ +\ {(J_H u)}^2\, [H_\P(Du)]^\bot H_{\P \P}(Du) \Big):D^2u\ = \ 0
\eeq
where $J_H$ is the \emph{$H$-Jacobian} of $u:\R^n \larrow \R^N$: 
\beq \label{eq13b}
J_H u \ := \ \left\{\begin{array}{l}
\det \left(H_\P(Du)^\top H_\P(Du) \right)^{\frac{1}{2}},\ \ \ n\leq N,\\
\det \left(H_\P(Du) H_\P(Du)^\top \right)^{\frac{1}{2}},\ \ \ n\geq N.
\end{array}
\right.
\eeq
When $H(\P)=\frac{1}{2}|\P|^2$, $J_H u$ is the ordinary Jacobian $Ju$ (see e.g.\ \cite{S}) and \eqref{eq13a} reduces to
\beq \label{eq14a}
\Gamma_\infty u\, :=\, Du\ot Du : D^2u\ +\ {(Ju)}^2\, [Du]^\bot \De u\ = \ 0.
\eeq
Systems \eqref{eq13a}, \eqref{eq14a} have continuous coefficients and their solutions are \emph{immersed or submersed varifolds}: $Ju$ is continuous and vanishes when $\rk(Du)< \min\{n,N\}$. The range $u(\R^n)$ may have lower-dimensional ``portions'' but by the Area-Coarea Formulas is $\mH^{\min\{n,N\}}$-a.e.\ equal to a sequence of manifolds immersed or submersed by $u$. However, by forcing the normal term to vanish at critical points there is a certain ``loss of information'' and uniqueness questions should be asked in the ``a.e.'' sense. Notwithstanding, $\Gamma_\infty$ and $\De_\infty$ are \emph{equivalent}  in the case of curves for $n=1$, $N\geq 1$ and trivially in the scalar case for $n\geq 1$, $N=1$, but generally do not coincide for maps $u:\R^n \larrow \R^N$ (Examples \ref{ex1}, \ref{ex3}).

In Section \ref{section5} we present a local structural decomposition property of regular \emph{Aronsson Maps}, that is, solutions $u\in C^2(\R^n)^N$ of \eqref{eq13}. They split in \emph{two phases}: for every point $x$, $\R^N$ decomposes to the direct sum of an horizontal subspace which is the range of the operator $H_P(Du(x))$ and a vertical subspace which is the nullspace of $H_P(Du(x))^\top$, both parameterized by $\R^n$. Actually, this structure is a vector bundle. In the case of $\De_\infty$, those subspaces are the tangent space $R(Du(x))$ and the normal space $N(Du^\top(x))$ of the image $u(\R^n)$ at $x$ respectively. We show that cross sections valued in the horizontal sub-bundle give rise to flows along whose trajectories the Hamilton Jacobi equation $H(Du)=c$ is satisfied. Cross sections valued in the vertical sub-bundle are coefficients of a \emph{linear} 1st order system coupled by a Hamilton Jacobi PDE which $u$ solves. For the case of $\De_\infty$, the PDE is linear as well. The geometric interpretation is that Aronsson maps generate a geometric structure like a sort of \emph{Ehressman connection in a fiber bundle}, which is ``vertically flat". However, the split is not of constant rank and interfaces of jumps in the rank appear. Hence the decomposition in not a global property, not even in the smooth case. It relates to the phase separation observed for Tight Maps in \cite{SS} but is different.

A basic difficulty in the vector case is the emergence of \emph{singular ``solutions''} which can not be rigorously justified.  In the scalar case, Viscosity Solutions have solved this problem. In Section \ref{section6} we construct $2$-dimensional $\infty$-Harmonic local $C^1$ diffeomorhisms and by following Katzourakis \cite{K} $C^1$ Aronsson Maps in arbitrary dimensions. Both classes are interpreted as everywhere solutions of the \emph{contracted version of \eqref{eq13}} (see \eqref{b9}, \eqref{c3}) which however is \emph{not equivalent to \eqref{eq13}}. By employing any nowhere differentiable auxiliary function $K$ in the explicit formulas \eqref{b8}, \eqref{eq17}, such ``solutions" may be nowhere twice differentiable and then their Hessians are genuine distributions not realizable by Radon measures, since such $K$'s are not in $BV(\R)$. 

In view of the fact that \eqref{eq13} is non-linear non-monotone and in non-divergence form, generally it does not possess classical, strong, weak, measure-valued, distributional or viscosity solutions. Hence, the standing ad hoc $C^2$ assumption for the solutions in all the recent works on the vector case makes their results of limited applicability. However, under the current state of the art in PDE theory there does not seem to be any variant of classical approaches allowing to tackle system like the $\infty$-Laplacian effectively.

\section{Derivation of the Tangential Aronsson System for Absolutely Minimizing Maps and Gradient Flows.} \label{section2}

\subsection{The Tangential Aronsson PDE System for Absolutely Minimizing Maps.} \label{subsection2.1} We begin by recalling Aronsson's minimality notion adapted to the vector case. The map $u\in W^{1,\infty}_{loc}(\R^n)^N$ is \emph{Absolutely Minimizing} for \eqref{eq1} when for all $\Om \subset \subset \R^n$ and all
$g \in W_0^{1,\infty}(\Om)^N$ we have
 \beq  \label{equ9}
E_\infty(u,\Om) \ \leq \ E_\infty(u+g,\Om).
 \eeq
This local minimality is essential and can not be deduced by minimality, since $\Om \mapsto E_\infty(u,\Om)$ is not a measure. In this subsection we prove the following

 \bt \label{th1}
  Let $H \in C^1\big((\R^{N}\ot \R^n) \by \R^N \by \R^n \big)$, $H \geq
 0$. If $u :\R^n \larrow \R^N$ is an Absolutely Minimizing Map for the supremal functional \eqref{eq1} and moreover $u\in C^2(\R^n)^N$, then $u$ is a solution to the tangential Aronsson PDE system \eqref{eq8}.
 \et

In particular, Theorem \ref{th1} extends a related result of Capogna and Raich in \cite{CR} for $H(\P)=|\P|^n/\textrm{det}(\P)$ and $n=N$ and provides a rigorous proof of a variant of Theorem \ref{th1} of Barron, Jensen and Wang in \cite{BJW2} which was merely sketched.

\BPT \ref{th1}. We follow an idea of R.\ Jensen from the scalar  $\De_\infty$ for $N=1$, appearing in \cite{J}. The idea is to derive \eqref{eq8} from \eqref{equ9} by inserting an appropriate energy comparison function to extract information and it turns out that a quadratic variation of $u$ suffices. Fix $x\in \R^n$, $\e,\de \in (0,1)$, $\xi \in \mS^{N-1}$ and consider the test function $g :\R^n \larrow \R^N$ given by
  \beq \label{equ10}
g(z)\ :=\ \frac{\de}{2}\big(\e^2 -|z-x|^2\big) \xi.
  \eeq
Obviously, $g \in W_0^{1,\infty}(\mB_\e (x))^N \cap C^2(\R^n)^N$. By setting $w :=  u +  g$, we calculate:
 \begin{align}
 w_\al(x)\ &=  \ u_\al(x)\ + \ \frac{\de \e^2}{2}\xi_\al, \label{equ13}\\
D_i w_\al(z)\ &=  \ D_i u_\al (z)\ - \ \de \xi_\al (z-x)_i,\\
D^2_{ij}w_\al (z)\ &=  \ D^2_{ij}u_\al (z)\ - \de \xi_\al \de_{ij}, \label{equ15}\\
 \|g\|_{C^0(\mB_\e(x))}\ &=\ \sup_{\mB_\e(x)}\,(g_\al\, g_\al)^\frac{1}{2} \  \leq \ 1, \label{equ15a}\\
  \|Dg \|_{C^0(\mB_\e(x))} \ &=\ \sup_{\mB_\e(x)}\, (D_i g_\al \, D_i g_\al)^\frac{1}{2} \ \leq \ 1,
  \end{align}
Since $u \in C^2(\R^n)^N$ and $H$ is $C^1$, the function
$H(Du,u,\_) :  z \mapsto H\big(Du(z),u(z),z\big)
$ is in $C^1(\R^n)$. By Taylor's theorem, there exists $\si \in
C^1(\R^n)$ with $\si(z-x)=o(1)$ as $z\rightarrow x$ such that
 \begin{align} \label{equ17}
H(Du,u,\_)(z)\ = & \ H(Du,u,\_ )(x) \ + \ D\big(H(Du,u,\_ )\big)(x)(z-x) \\
 & + \ \si(z-x)|z-x|.  \nonumber
 \end{align}
Again by Taylor's theorem for 
$H(D w,u,\_)\ : \ z \mapsto H\big(D w(z),u(z),z\big)$,
there exists $\tau \in C^1(\R^n)$ with $\tau(z-x)=o(1)$ as
$z\rightarrow x$ such that
 \begin{align} \label{equ17}
H(D w,u,\_)(z)\ =& \ H(D w,u,\_)(x) \ + \ D\big(H(D w,u,\_)\big)(x)(z-x)\\
 & + \ \tau(z-x)|z-x|. \nonumber
 \end{align}
Now we employ \eqref{equ17} to expand $H(D w,w,\_)$ as
\begin{align} \label{equ19}
  H\big(D w(z),w(z),(z)\big)
                 \ =& \ H\big(D w(x),u(x),(x)\big)  \nonumber\\
                 \ & + \ H_{P_{\al i}}\big(D w(x),u(x),(x)\big)\, D^2_{ij}w_\al
                 (x)\, (z-x)_j  \nonumber\\
                 \ & + \ H_{\eta_\al}\big(D w(x),u(x),(x)\big)\, D_j u_\al(x)\, (z-x)_j
                 \\
                 \ & +\ \tau(z-x)|z-x| \ + \ \big[H(D w,w,\_) - H(D w,u,\_)\big](z). \nonumber
\end{align}
Hence, by \eqref{equ19} and in view of \eqref{equ13} and \eqref{equ15}, we obtain
\begin{align} \label{equ20}
   H\big(D w(z),w(z),(z)\big)
    \ =& \ H\big(Du(x),u(x),(x)\big)   \nonumber\\
                 \ & + \ H_{P_{\al i}}\big(Du(x),u(x),(x)\big)\, D^2_{ij}u_\al
                 (x)\, (z-x)_j  \nonumber\\
                 \ & - \ \de \, \xi_\al \, H_{P_{\al
                 i}}\big(Du(x),u(x),(x)\big)\, (z-x)_i  \\
                 \ & + \ H_{\eta_\al}\big(Du(x),u(x),(x)\big)\,
                 D_j u_\al(x)\, (z-x)_j \nonumber\\
                 \ & +\ \tau(z-x)|z-x| \ + \ \big[H(D w,w,\_) - H(D w,u,\_)\big](z). \nonumber
\end{align}
Thus, \eqref{equ20} gives
\begin{align} \label{equ21}
  H(D w,w,\_)(z) \ =& \ H(Du,u,\_)(x) \ + \ D\big(H(Du,u,\_)\big)(x)\, (z-x)  \nonumber\\
                 \ & - \ \de \, \xi^\top H_{P}(Du,u,\_)(x)\, (z-x) \  +\ \tau(z-x)|z-x|  \\
                 \ &+ \ \big[H(D w,w,\_) - H(D w,u,\_)\big](z). \nonumber
\end{align}
By definition of $w$, we have the estimate
 \beq \label{equ22}
\sup_{\mB_\e(x)}\, |w - u| \ \leq\ \frac{\de \e^2}{2}.
 \eeq
Since by assumption $H$ is $C^1$, \eqref{equ22} implies the
estimate
 \beq \label{equ23}
\sup_{\mB_\e(x)}\, \Big|H(D w,w,\_) - H(D w,u,\_) \Big| \ \leq\ M\frac{\de
\e^2}{2},
 \eeq
with $M>0$ a constant independent of $\e,\de$ which can be taken to be
\beq \label{equ24}
M \ : =\ \sup \Big\{|H_\eta| \ \Big| \ \mB_R(0) \by \mB_R(0) \by
\mB_1(x), \ R:=1 + \|u\|_{C^1(\mB_1(x))}\Big\}. 
\eeq
Now we estimate energies. By \eqref{eq1} and \eqref{equ17}, we have
\begin{align} \label{equ25}
E_\infty(u,\mB_\e(x))\ =&\ \ \underset{\mB_\e(x)}{\ess\, \sup}\,  H(Du,u,\_ ) \nonumber\\
                        =&\ \ H(Du,u,\_ )(x) \ +\ \max_{\{|z-x|\leq
                        \e\}}\Big[D\big(H(Du,u,\_ )\big)(x)(z-x) \\
                        &\hspace{70pt} + \ \si(z-x)|z-x|\Big]
                        \nonumber
\end{align}
which gives
\begin{align} \label{equ26}
 E_\infty(u,\mB_\e(x))\ \ \geq & \ H(Du,u,\_ )(x) \ +\ \max_{\{|z-x|\leq
                        \e\}}\Big[D\big(H(Du,u,\_ )\big)(x)(z-x)\Big]\\
                        &\ - \ \max_{\{|z-x|\leq \e\}} |\si(z-x)|\,
                        |z-x|.\nonumber
\end{align}
The first maximum appearing in \eqref{equ26} is realized at
 \beq \label{equ27}
z_\e\ := \ x\ + \ \e\, \sgn\Big(D\big(H(Du,u,\_ )\big)(x)\Big)
 \eeq
while the second one is a quantity with decay $o(\e)$ as $\e\ri 0$.
Hence, \eqref{equ26} implies
\begin{align} \label{equ28}
 E_\infty(u,\mB_\e(x))\ \ \geq  \ H(Du,u,\_ )(x) \ +\ \e \, \big|D\big(H(Du,u,\_ )\big)(x)\big| \ + \   o(\e).
\end{align}
By \eqref{equ21} we have
\begin{align} \label{equ29}
E_\infty(w,\mB_\e(x))\ =&\ \ \underset{\mB_\e(x)}{\ess\, \sup}\, H(Dw,w,\_ ) \nonumber\\
                        =&\ \ H(Du,u,\_ )(x) \nonumber\\
                        &\ +\ \max_{\{|z-x|\leq
                        \e\}}\Big[\Big(D\big(H(Du,u,\_ )\big) - \de\,
                        \xi^\top
                        H_\P(Du,u,\_ )\Big)(x)(z-x) \\
                        &\hspace{40pt} + \ \tau(z-x)|z-x| \ +\ \big[H(D w,w,\_)
                                 - H(D w,u,\_)\big](z) \Big]. \nonumber
\end{align}
By employing estimate \eqref{equ23}, \eqref{equ29} implies
\begin{align} \label{equ30}
E_\infty(w,\mB_\e(x))\ \leq & \ H(Du,u,\_ )(x)\ + \ \max_{\{|z-x|\leq
                        \e\}} |\tau(z-x)|\, |z-x| \ + \ M\de \e^2 \nonumber\\
                        & +\ \max_{\{|z-x|\leq
                        \e\}}\Big[\Big(D\big(H(Du,u,\_ )\big) - \de\,
                        \xi^\top
                        H_\P(Du,u,\_ )\Big)(x)(z-x)\Big]. 
\end{align}
The second maximum appearing in \eqref{equ30} is realized at
 \beq \label{equ31}
z_{\e,\de,\xi}\ := \ x\ + \ \e\, \sgn\Big(D\big(H(Du,u,\_
)\big)(x) - \de\, \xi^\top H_\P(Du,u,\_ )(x)\Big),
 \eeq
while the first one is a quantity $o(\e)$ as $\e \ri 0$.
Hence, \eqref{equ30} gives
 \begin{align} \label{equ31}
E_\infty(w,\mB_\e(x))\ \leq & \ H(Du,u,\_ )(x) \ + \ o(\e)\ + \ M\de \e^2\nonumber\\
                        & +\ \e \, \Big|D\big(H(Du,u,\_
                  )\big)(x) - \de\, \xi^\top H_\P(Du,u,\_ )(x)\Big|.
\end{align}
By energy estimates \eqref{equ28}, \eqref{equ31} and Absolute
Minimality of $u$, we have 
 \begin{align} \label{equ32}
0\ \leq &\ E_\infty(w,\mB_\e(x))\ - \ E_\infty(u,\mB_\e(x)) \nonumber\\
                     \leq & \ H(Du,u,\_ )(x) \  + \ o(\e)\ + \ M\de \e^2       \ -\ \e \, \big|D\big(H(Du,u,\_ )\big)(x)\big| \\
                     & +\ \e \, \Big|D\big(H(Du,u,\_
                  )\big)(x) - \de\, \xi^\top H_\P(Du,u,\_ )(x)\Big| \ - \ H(Du,u,\_ )(x). \nonumber
\end{align}
By \eqref{equ32} we obtain
 \begin{align} \label{equ33}
0\ \leq & \ \Big|D\big(H(Du,u,\_ )\big)(x)
              - \de\, \xi^\top H_\P(Du,u,\_ )(x)\Big| \\
              & -\ \big|D\big(H(Du,u,\_ )\big)(x)\big|\
               + \ o(1)\ + \ M\de \e,  \nonumber
\end{align}
as $\e \ri 0$. Passage to the limit as $\e \ri 0$ gives
 \begin{align} \label{equ34}
0\ \leq & \ \Big|D\big(H(Du,u,\_ )\big)(x)
              - \de\, \xi^\top H_\P(Du,u,\_ )(x)\Big|\  -\ \big|D\big(H(Du,u,\_ )\big)(x)\big|.
\end{align}
If $D\big(H(Du,u,\_ )\big)(x)=0$, then by contracting derivatives in \eqref{eq8}, $u$ satisfies $\A_{\top , \infty} u = 0$. If on the other hand $D\big(H(Du,u,\_ )\big)(x) \neq 0$, the function
  \beq \label{equ35}
p\  \mapsto \ \big|D\big(H(Du,u,\_ )\big)(x)
              + p\, \big|  - \big|D\big(H(Du,u,\_ )\big)(x)\big|
 \eeq
is $C^1$ near the origin of $\R^n$.  By Taylor expansion of
function \eqref{equ35} at $0 \in \R^n$ applied at
$p = - \de\, \xi^\top H_\P(Du,u,\_ )(x)$, inequality \eqref{equ34} implies
 \beq \label{equ36}
0\ \leq  \ - \de\, \big(\xi^\top H_\P(Du,u,\_ )(x)\big)^\top \left(\frac{D\big(H(Du,u,\_
               )\big)(x)}{\big|D\big(H(Du,u,\_ )\big)(x)\big|}\right)\ +\ o(\de),
\eeq
as $\de \ri 0$. By \eqref{equ36} we deduce
 \begin{align} \label{equ37}
\xi^\top \Big(H_\P(Du,u,\_ )\big) D\big(H(Du,u,\_
)\big)(x)\Big) \ \leq \ o(1),
 \end{align}
as $\de \ri 0$. By passing to the limit, in view of \eqref{eq8} we obtain  $\xi^\top\big(\A_{\top , \infty} u(x)\big) \leq 0$. Since $\xi \in \mS^{N-1}$, $x \in \R^n$ are arbitrary, Absolutely Minimizing Maps $u:\R^n\larrow \R^N$ solve the tangential Aronsson PDE system and the theorem follows.
 \qed

\subsection{Gradient Flows for Tangentially $\infty$-Harmonic Maps.} Here we characterize regular solutions of the tangential $\infty$-Laplacian \ref{eq9} by employing appropriate gradient flows for a vector field depending on parameters (cf.\ \cite{CR}).  

\begin{proposition} \label{pr4}
Let $u \in C^2(\R^n)^N$. Then, $\De_{\top , \infty} u= 0$ on $\Om \subset \subset \R^n$ if and only if the flow map $\Phi : \R\by \Xi \larrow  \Om$ with $\Xi := \big\{ (x,\xi) \, | \, \xi^\top Du(x)\neq 0\big\}\sub \Om \by \mS^{N-1}$ of
\beq \left\{ \label{equ39}
\begin{array}{l}
\dfrac{\p}{\p t} \Phi(t,x,\xi) = \xi^\top  Du\big(\Phi(t,x,\xi)\big),\\
\ \ \ \Phi(0,x,\xi) = x,
\end{array}
\right.
\eeq
satisfies along trajectories
\beq \left\{ \label{equ38}
\begin{array}{c}
\big|Du\big(\Phi(t,x,\xi)\big)\big| = \big|Du(x)\big|, \ \ t \in \R, \ms\\
\, t\, \mapsto\, \xi^\top  u\big(\Phi(t,x,\xi)\big)\ \text{ is increasing}.
\end{array}
\right.
\eeq
\end{proposition}
\BPP \ref{pr4}. Follows directly from the following simple differential identities which can be verified by employing \eqref{equ39} and a straightforward calculation:
\beq
\frac{\p}{\p t}\left(\frac{1}{2}\big| Du\big(\Phi(t,x,\xi)\big)\big|^2\right)\ =\
\xi^\top \left(\De_{\top , \infty}  u\big(\Phi(t,x,\xi)\big)\right),
\eeq
\beq
\frac{\p}{\p t} \Big(\xi^\top  u\big(\Phi(t,x,\xi)\big)\Big)\ =\ \big| \xi^\top Du\big(\Phi(t,x,\xi)\big)\big|^2.
\eeq
We just note that if $\xi^\top Du(x)=0$, then in view of \eqref{eq9} we have $\De_{\top , \infty} u =0$. \qed

\ms \noi Unlike the case $N=1$ (see e.g. \cite{C}), images of trajectories are \emph{not} straight lines and we have to consider ODE systems with \emph{parameters}. It will be clear after Section \ref{section5} that an easy extension of Proposition \ref{pr4} holds for the general tangential Aronsson PDE system as well, if one assumes hypothesis \ref{b1} for $H$. Moreover, the next identities describe monotonicity and curvature of all projections $\eta^\top u$ along trajectories:
\begin{align}
\frac{\p}{\p t} \Big(\eta^\top  u\big(\Phi(t,x,\xi)\big)\Big)\ &= \ \Big[\xi \ot \eta :Du Du^\top\Big]\big(\Phi(t,x,\xi)\big),\\
\frac{{ \p }^2}{\p {t^2}} \Big( \eta^\top  u \big(\Phi(t,x,\xi)\big)\Big)\ &=\  \Big[D\big(\xi^\top u\big)^\top D\big(\xi \ot \eta :Du Du^\top \big)\Big] \big( \Phi(t,x,\xi)\big).
\end{align}
In particular, for $\xi=\eta$, 
\beq
\frac{{ \p }^2}{\p {t^2}} \Big( \xi^\top  u \big(\Phi(t,x,\xi)\big)\Big)\ =
\  2\, \De_\infty\big(\xi^\top u\big) \big( \Phi(t,x,\xi)\big).
\eeq
Hence, the projection $\xi^\top  u$ is affine along a trajectory if and only it is $\infty$-Harmonic thereon.

\section{The $\infty$-Laplacian with Respect to the Dual Operator Norm - Loss of Information in the Tangential Aronsson System.} \label{section3}

\subsection{Formal Derivation of the dual $\De_\infty$ for Tight Maps.}

Here we follow ideas from \cite{SS} and formally derive the dual version of the $\infty$-Laplacian which, as they prove, governs Tight Maps. Let $\|\_ \|_*$ be the operator norm on $(\R^N \ot \R^n)^*$ as in \eqref{eq10}. If $\P,\Q \in \R^N \ot \R^n$, by Danskin's theorem (\cite{D}), we have
\begin{align} \label{e1}
\frac{d}{dt}\Big|_{t=0}\Big(\frac{1}{2}\|\P + t\Q\|^2_*\Big) \ =\ \max_{\{\eta \in \mS^{N-1} : \|\P\|_* =|\eta^\top \P|\}} (\P \Q^\top) : \eta \ot \eta.
\end{align}
Let $u \in C^2(\R^n)^N$ and suppose $\|Du\|_*^2 = (DuDu^\top) :e\ot e$ for a unique $C^1$ unit vector field $e:\R^n \larrow \mS^{N-1}$, that is, that the maximum eigenvalue $\max \sigma (DuDu^\top)$ is a strict maximum of the spectrum. By setting $H(\P):=\frac{1}{2}\|\P\|_*^2$, \eqref{e1} implies
\beq \label{e2}
H_{\P_{\al i}}(Du) \ =\  e_\al e_\be D_i u_\be.
\eeq
The same choice of $H$ in \eqref{eq11} gives the following Euler-Lagrange PDE system
\beq \label{e3}
D_i\left(\Big(\frac{1}{2}\|Du\|_*^{2}\Big)^{p-1} e_\al e_\be D_i u_\be\right)\ =\ 0.
\eeq
Distributing derivatives, using \eqref{e2} and dividing by $(p-1)(\frac{1}{2}\|Du\|_*^{2})^{p-2}$, we obtain
\beq \label{e4}
(e_\be D_i u_\be)\, (e_\gamma D_i u_\gamma)\, (e_\de D^2_{ij} u_\de) e_\al \ =- \frac{\frac{1}{2}\|Du\|_*^{2}}{p-1}\Big[D_i (e_\be D_i u_\be)e_\al  + (D_i e_\al) \, (e_\gamma D_i u_\gamma) \Big]. 
\eeq
By expanding $(D_i e_\al)(e_\gamma D_i u_\gamma)$ to its projections on $\spn[e]$ and $(\spn[e])^\bot$ and using the identity
\beq \label{e5}
D_i \big(e_\al e_\gamma D_i u_\gamma\big)\ = \ e_\al \big((D_i e_\be)(D_i u_\be) + e_\be (D^2_{ii}u_\be ) \big)\ + \ (D_i e_\al)\, e_\be (D_i u_\be)
\eeq
we have
\begin{align} \label{e6}
De(e^\top Du)\ &= \ e\ot e \Big[ De(e^\top Du)\Big]\ +\ (I-e\ot e)\Big[ De(e^\top Du)\Big] \nonumber\\
& = \ \Big[ (De Du^\top):e\ot e \Big]e\ +\ e^\bot \Big[ \Div(e\ot e Du) -  (Du :De +e^\top \De u)e\Big] \\
& = \  \Big[ (De Du^\top):e\ot e \Big]e\ +\ e^\bot \Div(e\ot e\, Du). \nonumber
\end{align}
By \eqref{e4}, \eqref{e6} and perpendicularity we obtain by a further re-normalization that
\begin{align} \label{e7}
(e^\top Du)\ot (e^\top Du):&(e^\top D^2u) \, e\ \, +\, \ e^\bot \Div (e \ot e\, Du) \nonumber\\ 
&=\ -\frac{\frac{1}{2}|e^\top Du|^2}{p-1}\Big[(De Du^\top):e\ot e\, +\, \Div(e^\top Du)\Big].
\end{align}
Sending $p\ri \infty$ we obtain \eqref{eq12}. PDE system  \eqref{eq12} is equivalent to \eqref{eq5} obtained in \cite{SS}, but the respective form is different because we use the dual expression of the operator norm.

\subsection{The Loss of Normal Information in the Tangential Aronsson System.} Let now $H=H(Du)$ in \eqref{eq8}. If $\mS(n)$ denotes the symmetric matrices of $\R^n \ot \R^n$, we introduce the \emph{Contraction} operation
\beq  \label{e8} 
(\R^N \ot \R^n) \ot (\R^N \ot \R^n) \by \R^N \ot \mS(n) \larrow \R^N \ :\ 
(\textbf{A} :\X)_\al\; :=\; \textbf{A}_{\al i \be j} \X_{\be ij}.
\eeq
This operation extends the usual trace inner product to higher order tensors. In view of \eqref{e8}, we may contract \eqref{eq8} to
\beq \label{e9}
\A_{\top , \infty} u \ = \ H_\P (Du) \ot H_\P(Du):D^2u\ =\ 0.
\eeq
Using \eqref{e2} into \eqref{e9}, we readily obtain 

\begin{proposition} Let $H(\P)=\frac{1}{2}\|\P\|^2_*$ and suppose $u \in C^2(\R^n)^N$ and that there exists $e\in C^1(\R^n)^N$ with $|e|\equiv 1$ such that $\|Du\|_*^2=(DuDu^\top):e\ot e$. Then,
\beq \label{e10}
\A_{\top , \infty} u  \ = \ e\ot e\, \big( \De_{\infty_{T^*}} u\big). 
\eeq
Thus, the tangential Aronsson system \eqref{e9} which governs Absolutely Minimizing Maps equals the projection on $\spn[e]$ of the $\infty$-Laplacian \eqref{eq12} which governs Tight Maps.
\end{proposition}
This loss of information along the hyperplane $(\spn[e])^\bot$ renders \eqref{e9} defective from the viewpoint of uniqueness for the Dirichlet problem and sufficiency for the variational problem:
\begin{example} \label{ex1}(cf. \cite{OTW}, \cite{WO}) The Dirichlet problem for the tangential $\De_{\top , \infty} $
\[
\left\{\begin{array}{l}
\De_{\top , \infty} u\ = \ 0,\ \ \ u:(0,1)\larrow \R^N,\\
u(0)\, =\, u(1)\, =\, 0, \ \ N>1,
\end{array}
\right.
\]
admits infinitely many smooth solutions. Actually, the class of tangentially $\infty$-Harmonic curves coincides with the class of constant speed curves. To see this, fix $u\in C^2(\R)^N$. Then, by \eqref{eq9} we have
 \beq
\De_{\top , \infty} u\ =\ (\dot{u}\ot \dot{u})\ddot{u} \ =\ \Big(\frac{1}{2}|\dot{u}|^2\Big)^. \dot{u}.
 \eeq
If $|\dot{u}|^2$ is constant, then we have $\De_{\top , \infty} u =0$. Conversely, if $\De_{\top , \infty} u=0$, then we have that $|\dot{u}|^2$ is constant on each connected component of the open set $\{\dot{u} \neq0\}$ but the same holds also on its complement $\{\dot{u}=0\}$. Moreover, $|\dot{u}|^2 \in C^1(\R)$. Consequently, $|\dot{u}|^2$ is constant since it can not have discontinuity jumps.  
\end{example}

\begin{example} \label{ex2}
\eqref{e9} fails to characterize Absolutely Minimizing Maps, even in the case of the tangential $\infty$-Laplacian $\De_{\top , \infty} $ for curves on the plane: let $n=1$, $N=2$, $H : \R^2 \larrow \R$ given by $H(\P)=|\P|^2$, $\Om = (0,2\pi)$, and
\beq
u(x)=(\cos x, \sin x)^\top\ ,\ \ \  v(x)=(1,0)^\top.
\eeq
Then, $u=v$ on $\p \Om$ and $\De_{\top , \infty} u=\De_{\top , \infty} v=0$ on $\Om$ but 
\begin{align}
E_\infty(u,\Om)\ &=\ 1 \nonumber\\
&>\ 0\ =\ E_\infty(v,\Om).  \nonumber
\end{align}
Hence $u$ is not a minimizer of $E_\infty$ over $\Om$, since its energy can improve further.
\end{example}

\section{Recovering the Lost Information - The Full Aronsson System.} \label{section4}

\subsection{Derivation of the Aronsson PDE system in the $L^p$ limit.} \label{subsection1} Here we restrict attention to $H=H(\P) \in C^2(\R^N \ot \R^n)$ and formally derive the complete Aronsson PDE system. The discussion applies to general $H=H(\P,\eta,x)$ equally well. In view of \eqref{e8}, the Euler Lagrange PDE system of \eqref{eq11} after expansion and division by $p(p-1){H(Du)^{p-2}}$ is
\beq \label{a1}
H_\P (Du) \ot H_\P(Du):D^2u \ + \ \frac{H(Du)}{p-1}H_{\P\P}(Du):D^2u\ = \ 0.
\eeq
Let $\textrm{N}\big(H_\P (Du(x))^\top\big)$ denote the null space of $H_\P (Du(x))^\top :\R^N\larrow \R^n$ and $\textrm{R}\big(H_\P (Du(x))\big)$ the range of $H_\P (Du(x)) :\R^n \larrow \R^N$. We define the projections
\begin{align} 
[ H_\P (Du) ]^\bot \ & := \ \textrm{Proj}_{\, \textrm{N}(H_\P (Du)^\top)} \label{a2a}\\
[ H_\P (Du) ]^\top \ & := \ \textrm{Proj}_{\, \textrm{R}(H_\P (Du))} .  \label{a2b}
\end{align}
When $N=1$ and $H_\P (Du)\neq0$, we have $[H_\P (Du)]^\bot \equiv 0$ and $[H_\P (Du)]^\top \equiv 1$. Generally, for $N\geq 1$ we have $[ H_\P (Du) ]^\bot + [ H_\P (Du) ]^\top =I$. By expanding the term $H_{\P\P}(Du):D^2 u$ with respect to these projections and using that $[ H_\P (Du) ]^\top H_\P (Du)= H_\P (Du)$, \eqref{a1} gives

\begin{align} \label{a3}
H_\P (Du) \ot H_\P(Du):D^2u \ + \ & \frac{H(Du)}{p-1}[ H_\P (Du) ]^\top H_{\P\P}(Du):D^2u \nonumber\\
 &=\ - \ \frac{H(Du)}{p-1}[ H_\P (Du) ]^\bot H_{\P\P}(Du):D^2u.
\end{align}
By perpendicularity, both sides of \eqref{a3} vanish. We re-normalize the right hand side by multiplying by $p-1$ and rearrange, to obtain 
\begin{align} \label{a4}
\Big(H_\P (Du) \ot H_\P(Du)\ +\  & H(Du)[ H_\P (Du) ]^\bot H_{\P\P}(Du)\Big):D^2u \nonumber\\
 &=\ - \ \frac{H(Du)}{p-1}[ H_\P (Du) ]^\top H_{\P\P}(Du):D^2u.
\end{align}
In the limit as $p\ri \infty$ we obtain the full Aronsson PDE system \eqref{eq13}. The choice $H(\P)=\frac{1}{2}|\P|^2$ gives the full $\infty$-Laplacian \eqref{eq14} with respect to the Euclidean norm. In the scalar case of $N=1$ the normal term vanishes identically and scalar functions solving  Aronsson's PDE \eqref{eq13} coincide with functions solving Aronsson's  tangential PDE \eqref{eq8}.
\begin{example} \label{ex3}
The additional term remedies the non-uniqueness phenomena of Example \ref{ex1}. Let $n=1$, $N\geq 1$ and $u\in C^2(\R)^N$ non-constant. For $n=1$, the ODE systems \eqref{eq14} and \eqref{eq14a} are equivalent since solutions have no critical points. So, $\De_\infty$ and $\Gamma_\infty$ are \emph{equivalent}. Hence, we have
\beq
\begin{array}{l}
\De_\infty u\ = \ {\dot{u}}\ot {\dot{u}}:{\ddot{u}}\ + \ |\dot{u}|^2\,  [\dot{u}]^\bot \ddot{u} \\
\hspace{30pt}=\ ({\dot{u}}\ot {\dot{u}}){\ddot{u}}\ +\ |\dot{u}|^2 \Big(I\, - \, \dfrac{\dot{u}}{|\dot{u}|}\ot  \dfrac{\dot{u}}{|\dot{u}|} \Big)\, \ddot{u}\\
\hspace{30pt}=\ |\dot{u}|^2 \,\ddot{u}. 
\end{array} 
\eeq
Consequently, $\infty$-Harmonic curves are affine and their graphs are straight lines.
\end{example}
Inspired from \cite{SS}, we propose the following \emph{minimality notion for general $L^\infty$ variational problems}: let $H \in C^2(\R^N \ot \R^n)$, $H\geq H(0)=0$ rank-1 convex. For $E_\infty(u,\Om)=\ess \, \sup_\Om H(Du)$, we set $E_\infty u(x) := \lim_{r\ri0}E_\infty(u,\mB_r(x))$. Then, the map $u :\Om \subset \subset \R^n \larrow \R^N$ is \textbf{Stark} if $u\in W^{1,\infty}(\Om)^N$ and for all $v\in W_u^{1,\infty}(\Om)^N$,
\beq
\sup_{\{E_\infty u > E_\infty v\}} E_\infty u \ \leq \ \sup_{\{E_\infty v > E_\infty u\}} E_\infty v.
\eeq
We conjecture that Starkness characterizes Aronsson Maps, i.e. solutions of \eqref{eq13}.

\subsection{Smooth $\infty$-Harmonic Maps with interfaces of rank discontinuities (cf.\ \cite{OTW},\cite{WO})} \label{subsection2} The normal coefficient $|Du|^2[Du]^\bot$ of $\De u$ generally is discontinuous, no matter how smooth the solution might be since rank jumps of $Du$ occur even for $C^\infty$ solutions. Let $N\geq 2$ and $f,g:\R \larrow \R^N$ smooth curves and define $u:\R^2 \larrow \R^N$ by
\beq \label{a5}
u(x,y)\ :=\ f(x)\ + \ g(y).
\eeq
Then, $Du(x,y)=(f'(x),g'(y))$ which is a tensor of $\R^N \ot \R^2$ and $\textrm{N}\big(Du(x,y)^\top\big)=\big(\spn[\{f'(x),g'(y)\}]\big)^\bot$. By a simple calculation,
\begin{align} \label{a6}
\De_\infty u(x,y)\ = \ & \left(f'(x) \ot  f'(x)\right)f''(x)\ \ +\  \  \big(g'(y) \ot g'(y)\big)g''(y)\\ 
& + \ \ \big(|f'(x)|^2 +|g'(y)|^2\big)[\big(f'(x),g'(y)\big)]^\bot\Big(f''(x)\, +\, g''(y)\Big)  \nonumber.
\end{align}
Suppose further that $f$, $g$ are parametric unit speed curves. The normal coefficient is discontinuous on the interface $\mathcal{S}$ whereon $f'(x)$ becomes co-linear to $g'(y)$: for $[f'(x)]^\bot = I- f'(x)\ot f'(x)$ and
\beq
\mathcal{S}\ :=\ \big\{(x,y)\ :\ |f'(x)^\top g'(y)|=1\big\}
\eeq
we have
\beq
[\big(f'(x),g'(y)\big)]^\bot = 
 \left\{
\begin{array}{l}
I - f'(x) \ot f'(x)\, -\, \dfrac{[f'(x)]^\bot g'(y)}{\big|[f'(x)]^\bot g'(y)\big|} \ot \dfrac{[f'(x)]^\bot g'(y)}{\big|[f'(x)]^\bot g'(y)\big|}  ,\ms \\
\hspace{196pt}(x,y)\in \R^2 \set \mathcal{S}, \ms\\
I - f'(x) \ot f'(x), \hspace{140pt} (x,y)\in \mathcal{S}.
\end{array}
\right. 
\eeq
In particular, let $N=2$ and $g=\pm f$. Then, by \eqref{a6} we have 
\begin{align} \label{aa6}
\De_\infty u(x,y)\ = \  2[\big(f'(x),\pm f'(y)\big)]^\bot\Big(f''(x)\, \pm\, f''(y)\Big) .
\end{align}

On the diagonal we have $Du(x,x)=f'(x)\ot (e_1\pm e_2)$, the rank decreases from $2$ to $1$ but $Du(x,x)$ does not vanish and $[Du(x,y)]^\bot$ in discontinuous on $\{x=y\}\sub \mathcal{S}$. What happens is that the tangent space of the graph of $u$ is spanned by $2$ smoothly moving vectors, which become colinear along the interface.
\[
\begin{array}{c}
\includegraphics[scale=0.15]{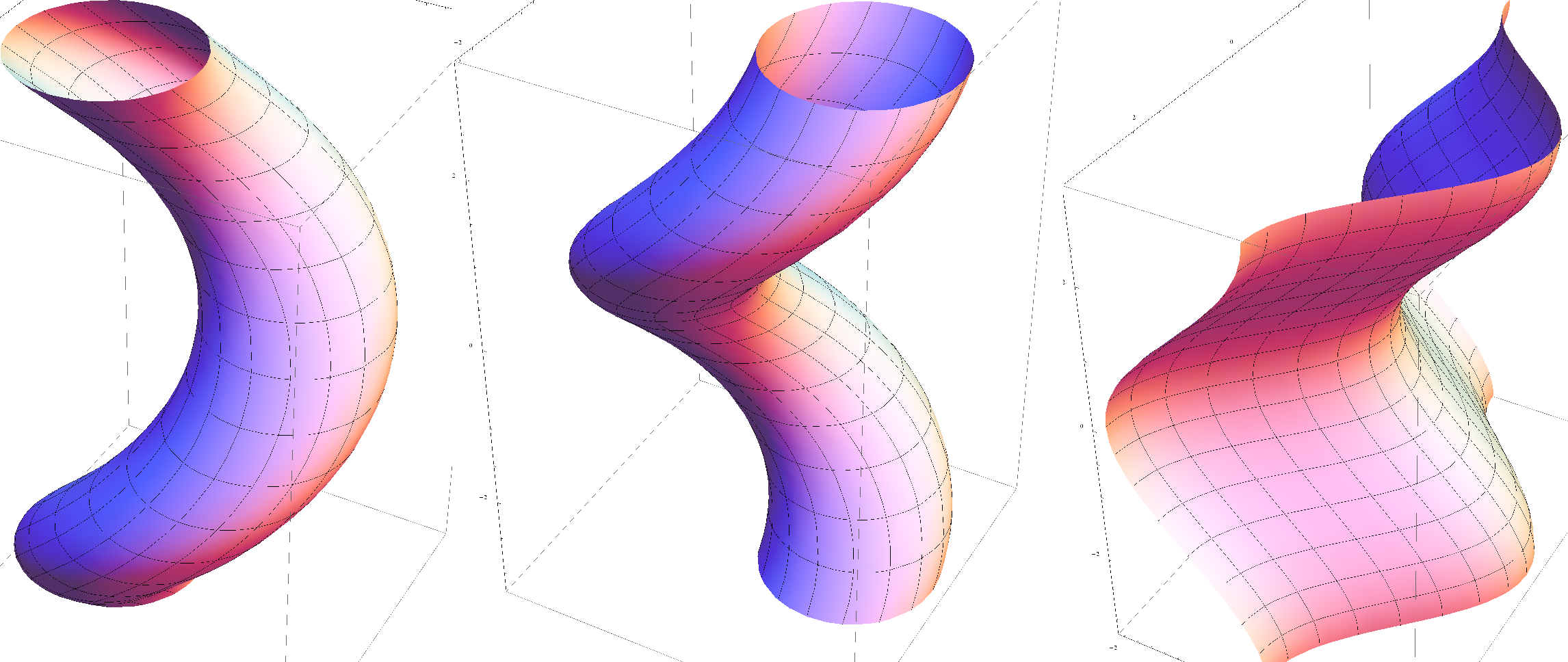}\\
\underset{\text{3 projections on $\R^3$ of the graph of $u(x,y)=e^{ix}-e^{iy}$, its range on $\R^2$ and its ``covering sheets''.}}
{\includegraphics[scale=0.18]{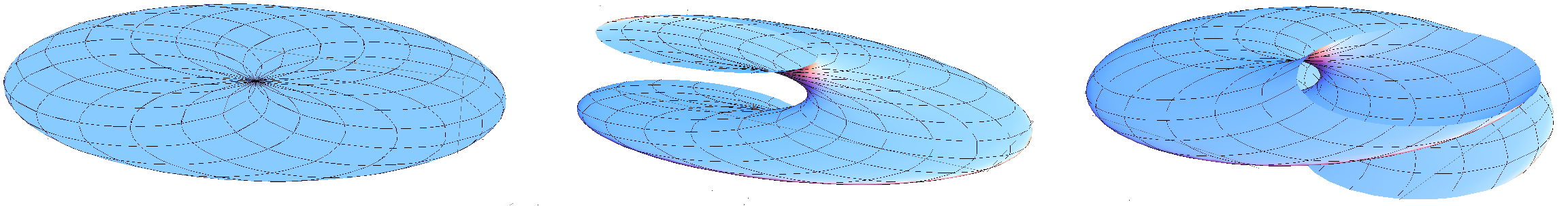}}
\end{array}
\]
 If $g=f$, the normal term is \emph{not} annihilated so $u$ is \emph{not $\infty$-Harmonic}   (cf.\ \cite{OTW}).  This is compatible with Aronsson's scalar example (\cite{A6}): $(x,y)\mapsto|x|^{\frac{4}{3}}-|y|^{\frac{4}{3}}$ is $\infty$-Harmonic, but $(x,y)\mapsto|x|^{\frac{4}{3}}+|y|^{\frac{4}{3}}$ is \emph{not}. 

The choice $f(t)=e^{it}$ gives the $\infty$-Harmonic map $u(x,y)=e^{ix}-e^{iy}$, $u\, : \, \Om=\{|x\pm y|<\pi\} \sub \R^2 \larrow \R^2$. Since $f'(x)^\top g'(y)=\cos(x-y)$, $u$ solves \eqref{a6} on the rhombus $\Om$ with $\mathcal{S}=\{x=y\}\cap \Om$. Since $\De_\infty u(0,\pi)=-4e_1 \neq 0$, $u$ is not a solution outside of the strip $\{|x-y|<\pi\}$. $\infty$-Harmonic Maps with discontinuous rank in their gradients and hence discontinuities in the dimension of their tangent space are varifolds, rather than manifolds. Since in the measure-theoretic approach lower dimensional portions of $u(\R^n)$ are neglected and they seem to matter, it is perhaps better to be viewed as Cell Complexes, a structure well studied in Algebraic Topology.

The varifold modification \eqref{eq14a} of \eqref{eq14} is legitimate and rather natural in order to produce $C^0$ coefficients, but not alter $\infty$-Harmonic maps much. On zeros of $Du$ and off interfaces of rank discontinuities, both systems are equivalent. On the interfaces but off zeros, multiply the normal PDE system $|Du|^2[Du]^\bot\De u=0$ by $({Ju}/{|Du|})^2$ and this forces discontinuities to vanish. Such phenomena do not appear for $n=1$ or $N=1$, since then $Ju$ equals either $|\dot{u}|$ or $|Du|$ respectively.

\section{Phase Decomposition of Aronsson Maps -  Gradient Flows Vs Linear 1st Order Systems Coupled by Hamilton Jacobi PDE.} \label{section5}

In this section under a mild radial monotonicity assumption on the Hamiltonian we present a structural decomposition property of Aronsson Maps. There is a local splitting to an horizontal vector bundle with fibers  $R\big(H_P(Du(x))\big)$ and a vertical vector bundle with fibers  $N\big(H_P(Du(x))^\top\big)$ whose direct sum is the trivial bundle $\R^N \oplus \R^n \larrow \R^n$. Horizontal sections generate flows such that along integral curves $u$ solves the PDE $H(Du)=c$, while vertical sections furnish coefficients of a \emph{linear} 1st order systems coupled by a Hamilton Jacobi PDE which $u$ solves. However, in view of Subsection \ref{subsection2} \emph{the bundles are not of constant rank} and \emph{interfaces of discontinuities in the dimension of the fibers appear}. Hence, even in the smooth case this local structure generally does not extend to a global one.

\subsection{The Phase Decomposition Theorem.} Suppose that $H \in C^2(\R^N \ot \R^n)$ is \emph{radially rank-1 monotone increasing}, in the sense that there exists a positive increasing $\omega \in C^0(0,\infty)$ with $\omega(0^+)=0$ such that
\beq \label{b1}
H_\P (\Q)\Q^\top : \xi \ot \xi \ \geq \ \omega(\big| \xi^\top \Q \big|),
\eeq
for all $\Q \in \R^N \ot \R^n$ and $\xi \in \mS^{N-1}$. When $H(\P)=\frac{1}{2}|\P|^2$, \eqref{b1} holds for $\omega(t)=t^2$.

\bt \label{th2} Let $H \in C^2(\R^N \ot \R^n)$ satisfy assumption \eqref{b1} and suppose $u \in C^2(\R^n)^N$ is an Aronsson Map solving \eqref{eq13}. Then, for any open set $\Om \sub \R^n$ whereon $H_\P(Du)$ has constant rank and any $C^1$ unit vector field $e :\R^n \larrow \mS^{N-1}$ there exists an orthogonal split $e=h+v$ to horizontal and vertical component and an open cover of $\Om$ by a pair $\{\Om^\top,\Om^\bot\}$ given by
\begin{align}
\Om^\top = \{x\in \Om \ :\ h(x)\neq 0\}\ , \ \ \Om^\bot = \{x\in \Om \ :\ v(x)\neq 0\},
\end{align}
such that the flow map $\Phi :\R \by \Om^\top \larrow \Om^\top$ of
\beq \label{b2}
\left\{
\begin{array}{l}
\dfrac{\p}{\p t} \Phi(t,x)\ = \ h\big(\Phi(t,x)\big), \ms \\
\ \ \ \Phi(0,x)\ =\ x,  
\end{array}
\right.
\eeq
satisfies along trajectories
\beq \begin{array}{l}
 H\left(Du\big(\Phi(t,x)\big)\right)\ = \ H(Du(x)), \label{b3} \ms\\
  h\big(\Phi(t,x)\big)^\top \dfrac{\p}{\p t} \big(u\big(\Phi(t,x)\big)\big) \ \geq\ 0,
\end{array}
\eeq
and also $u$ solves on $\Om^\bot$ the  PDE system
\beq \label{b5}
\left\{\begin{array}{c}
v^\top Du \, =\, 0,\\
Dv:H_\P(Du)\, =\, 0.
\end{array}\right.
\eeq
\et
\noi We note that in the case of quadratic Hamiltonians, \eqref{b5} becomes a \emph{linear} system. In particular, for $\De_\infty$ we obtain 
\beq
v^\top Du\ =\ 0\ ,\ \ Dv: Du\ =\ 0
\eeq
and the respective vector bundles are the tangent and the normal respectively to the image of $u$ respectively. We remark that the coefficients \emph{depend on $Du$}, but we have the \emph{form} of a linear PDE system.

\BPT \ref{th2}. Suppose $\A_\infty u=0$, fix $e:\R^n \larrow \mS^{N-1}$ in $C^1$ and define 
\begin{align}
h\, :=\, [H_\P(Du)]^\top e\ , \ \ \ \ v\, :=\, [H_\P(Du)]^\bot e.  \label{b7}
\end{align}
By \eqref{eq13}, we have
\begin{align}
h^\top\big(H_\P(Du) \ot H_\P(Du) :D^2u\big)\ =\ h^\top \A_\infty u \ = \ 0,\label{ab1}\\
v^\top\big([H_\P(Du)]^\bot H_{\P\P}(Du) :D^2u\big)\ =\ v^\top \A_\infty u \ = \ 0.\label{ab2}
\end{align}
We use \eqref{b2}, \eqref{ab1} and that $h$ is horizontal and calculate
\begin{align}
\frac{\p }{\p t}H\big(Du(\Phi(t,x))\big) &= H_{\P_{\al i}}\big(Du(\Phi(t,x))\big)\, D^2_{ij}u_\al (\Phi(t,x))\, \frac{\p }{\p t}\Phi_j(t,x)\nonumber\\
&= \Big[h_\be\, H_{\P_{\be j}}(Du)\,H_{\P_{\al i}}(Du)\, D^2_{ij}u_\al \Big](\Phi(t,x))\\
&= h^\top \big(\A_\infty u(\Phi(t,x))\big). \nonumber
\end{align}
We now use \eqref{b1}, \eqref{b2} and calculate
\begin{align}
h(\Phi(t,x))^\top \frac{\p }{\p t}\big(u(\Phi(t,x)\big) &= h(\Phi(t,x))^\top Du(\Phi(t,x))\, \frac{\p }{\p t}\Phi(t,x) \nonumber\\
&= h(\Phi(t,x))^\top Du(\Phi(t,x))\; h(\Phi(t,x))^\top H_\P \big(Du(\Phi(t,x))\big) \nonumber\\
&= H_\P \big(Du(\Phi(t,x))\big)Du(\Phi(t,x))^\top : h(\Phi(t,x)) \ot h(\Phi(t,x)) \\
&\geq \omega(\big|h(\Phi(t,x))^\top Du(\Phi(t,x)) \big|).  \nonumber
\end{align}
By \eqref{b7} we have that $v$ is the projection of $e$ on the null space of $H_\P(Du)$. Hence, $v^\top H_\P(Du)=0$ on $\Om^\bot$ and consequently by \eqref{b1}
\begin{align}
0\ &= \ \big( v^\top H_\P(Du)\big)\, \big(v^\top Du\big) \nonumber\\
&= \ H_\P(Du)Du^\top :v \ot v\\
&\geq \ \omega(\big|v^\top Du\big|). \nonumber
\end{align}
Thus, $v^\top Du =0$. By employing \eqref{ab2} and the differential identity
\beq
\Div\big(v^\top H_\P(Du)\big)\ =\ Dv:H_\P(Du) \ + \ v^\top\big(H_{\P\P}(Du):D^2 u\big) 
\eeq
we complete the derivation of \eqref{b5} and the theorem follows. \qed

\section{Singular $\infty$-Harmonic Local Diffeomorhisms and Aronsson Maps} \label{section6}

\subsection{$\infty$-Harmonic $C^1$ Immersions.} Let $K\in C^0(\R)$ and define $u : \R^2 \larrow \R^2$ by
\beq \label{b8}
u(x,y)\ :=\, \int_0^x e^{i K(t)}dt \ + \ i\int_0^y e^{i K(s)} ds.
\eeq 
Formula \eqref{b8} defines a $2$-dimensional $\infty$-Harmonic Map, which is singular unless $K\in C^1(\R)$.
\begin{proposition} \label{p2} Suppose $\|K\|_{C^0(\R)}< \frac{\pi}{4}$ and let $u$ be given by \eqref{b8}. Then, $u$ is a $C^1(\R^2)^2$-local diffeomorphism and everywhere solution of the PDE system
\beq \label{b9}
Du\, D\left(\frac{1}{2}|Du|^2\right)\ + \ |Du|^2[Du]^\bot \Div\, (Du)\ = \ 0.
\eeq
\end{proposition}
\BPP \ref{p2}. Let $f(x)$, $g(y)$ denote the two summands comprising \eqref{b8}. Trivially, $f,\, g$ are parametric unit speed curves. Since $|K| < \frac{\pi}{4}$, we obtain $|K(x)-K(y)| <\frac{\pi}{2}$ and consequently $f'(x)$ and $g'(y)$ are nowhere co-linear:
\beq
\big|f'(x)^\top g'(y)\big|\ = \ \big|\sin\big(K(x)-K(y)\big)\big| \ < \ 1,\ \ \ (x,y)\in \R^2.
\eeq
Thus, $\rk(Du(x,y))=2$, $[Du(x,y)]^\bot \equiv 0$ in $\R^2$ and the coefficient of the second summand in \eqref{b9} vanishes identically. By contracting derivatives in \eqref{a6} we have
\begin{align} \label{eq}
\De_\infty u(x,y)\ = \ & \left(\frac{1}{2}|f'|^2\right)'(x)\, f'(x)\ +  \ \left(\frac{1}{2}|g'|^2\right)'(y)\, g'(y)\ = \ 0.
\end{align}
Hence, $u$ is an everywhere solution of the contracted PDE \eqref{b9}.\qed

\ms \noi It is not difficult to see that Proposition \ref{p2} readily extends to higher dimensions.

\subsection{Singular Aronsson Maps.}

We follow Katzourakis \cite{K} and construct singular solutions to \eqref{eq13} under the assumption that some level set of $H$ contains a straight line segment of rank-1 matrices and has constant gradient $H_\P$ thereon. Given $a,b \in \R^n$, $\eta\in \R^N$ and $K\in C^0(\R)$, we define $u:\R^n \larrow \R^N$ by
\beq \label{eq17} 
u(x)\ :=\ \Bigg(x^\top \Big(\frac{b+a}{2}\Big) \ +\ \int_0^{x^\top \big(\frac{b-a}{2}\big)} K(t)dt \Bigg) \, \eta.
\eeq

\begin{proposition} \label{p1}
Let $H\in C^2(\R^N\ot \R^n)$ and suppose that there exist $a,b \in \R^n$, $\eta\in \R^N$, $c\in\R$ and $\textbf{C} \in \R^N\ot \R^n$ such that
\beq  \label{c1}
[\eta \ot a,\, \eta\ot b] \ \sub \ \{H=c\}\cap\{H_\P =\textbf{C}\}.
\eeq
Let $u$ be given by \eqref{eq17} and suppose $\|K\|_{C^0(\R)}<1$. Then, $u$ is an everywhere solution in $C^1(\R^n)^N$ of the PDE system
\beq \label{c3}
H_\P(Du)\, D\big(H(Du)\big)\ + \ H(Du)[H_\P(Du)]^\bot \Div \, \big(H_\P(Du)\big)\ = \ 0.
\eeq
\end{proposition}

\BPP \ref{p1}. We recall that $\|K\|_{C^0(\R)}<1$ and calculate
\begin{align} \label{c2}
Du(x)\ &=\ \eta \ot \frac{a+b}{2}\ + \ K\left(x^\top\Big(\frac{b-a}{2}\Big)\right)\eta \ot  \frac{b-a}{2} \nonumber\\
&=\  \left(\frac{1}{2} - \frac{1}{2}K\left(x^\top\Big(\frac{b-a}{2}\Big)\right)\right)   \eta \ot a\\ 
&\ \ \ \ \ + \ \left[1-\left(\frac{1}{2} - \frac{1}{2} K\left(x^\top\Big(\frac{b-a}{2}\Big)\right)\right)\right]\eta \ot b.  \nonumber
\end{align}
Hence, \eqref{c2} implies $Du(\R^n) \sub [\eta \ot a,\, \eta\ot b] $. For $u$ as in \eqref{eq17}, we have $H(Du)\equiv c$ and $H_\P(Du) \equiv \textbf{C}$ on $\R^n$ as a result of \eqref{c1}. Hence $u$ solves the system \eqref{c3}. 
\qed 
\begin{remark} The problem is evident. For $C^2$ solutions, we can contract \eqref{eq13} to \eqref{c3}, but they are \emph{not} equivalent. By choosing as $K$ in \eqref{b8} and \eqref{eq17} any nowhere differentiable function, $D^2u$ exists nowhere on $\R^n$, is not a Radon measure and can be seen only as a first order distribution. In the scalar case $N=1$, \eqref{b8} degenerates to an affine function and for the respective to \eqref{eq13} solutions to the Aronsson PDE it is proved in \cite{K} that they are viscosity solutions. However, if $N>1$, \eqref{eq13} is a quasilinear non-monotone PDE system in non-divergence form and by the previous facts generally does not possess classical, strong, weak, measure-valued, distributional or viscosity solutions.
\end{remark}

\ms

\noi \textbf{Acknowledgement.} I thank L.\ C.\ Evans, R.\
Jensen, Ch. Wang and J. Manfredi for their interest. Subection \ref{subsection2.1} of this work was written when the Author was a doctoral student at the Department of Mathematics, University of Athens, Greece.

\bibliographystyle{amsplain}

\end{document}